\documentclass[twocolumn]{article}
\usepackage{amsmath}
\usepackage{amsfonts}
\usepackage{amsthm}

\usepackage{graphicx,geometry}
\title{Solving the continuous nonlinear resource
allocation problem with an interior point method}

\author{
Stephen E. Wright\thanks{
  Department of Statistics,
  Miami University,
  Oxford, OH 45056.
  Email: {\texttt wrightse@muohio.edu}
  }
\ and
James J.~Rohal\thanks{
  Department of Mathematics, North Carolina State
   University, Raleigh, North Carolina 27695.
   Email: {\texttt jjrohal@ncsu.edu}}
}

\begin{document}

% comment this out for Elsevier
\maketitle

% needed for Elsevier
%\begin{frontmatter}

\begin{abstract}
Resource allocation problems are usually solved with specialized methods
exploiting their general sparsity and problem-specific algebraic structure.
We show that the sparsity structure alone yields a closed-form Newton search
direction for the generic primal-dual interior point method. Computational
tests show that the interior point method consistently outperforms the best
specialized methods when no additional algebraic structure is available.
\end{abstract}

% use these for Elsevier
%\begin{keyword}
%convex programming \sep interior point methods \sep continuous knapsack
%\end{keyword}

% needed for Elsevier
%\end{frontmatter}

\section{Introduction} \label{sec:introduction}

\newcommand{\diag}{\mathop{\mathrm{diag}}}
\newcommand{\argmin}{\mathop{\mathrm{argmin}}}

We consider the resource allocation problem in the form
\begin{align}
  \text{minimize } &f(x):=\sum_{i=1}^nf_i(x_i)
     \text{ over all } x\label{obj}\\
  \text{subject to }&g(x):=\sum_{i=1}^ng_i(x_i)=b,\label{gcon}\\
  & l\leq x\leq u.\label{intervals}
\end{align}
Here $x$, $l$, and $u$ are $n$-vectors of real numbers, $b$ is a real scalar,
and the functions $f_i$ and $g_i$ are convex and twice differentiable on an
open set containing the interval $[l_i,u_i]$. Inequalities of vectors are
interpreted coordinate-wise.

The recent survey paper of Patriksson \cite{Patriksson} shows that such
problems have a long history and diverse applications. The contexts in which
the problem appears often demand that it be solved very quickly, even in high
dimensions. Consequently, researchers long ago moved beyond general-purpose
nonlinear programming procedures and focused on exploiting the special
structure of the optimality conditions for the problem.  As noted by
Patriksson, two frameworks have emerged as the most competitive for solving
resource allocation problems: the \emph{pegging} or \emph{variable-fixing}
methods and the \emph{breakpoint-search} methods. Patriksson also observes
that computational studies in the literature have generally indicated that
pegging is superior to breakpoint search when certain subproblems (see
\S\ref{section:descr}) common to both methods are easily solved, whereas
breakpoint search is faster otherwise.  Moreover, numerical comparisons of
either method with general-purpose solvers are essentially nonexistent in the
literature.

Here we present evidence that a primal-dual interior point method outperforms
breakpoint search on problems for which the latter is traditionally
considered the best possible choice, namely, when its subproblems do not
admit closed-form solutions and must be solved numerically. We show that the
special structure of (\ref{obj})--(\ref{intervals}) allows for a closed-form
solution of the linear system defining the search directions and we present
computational results showing the method's superiority. This addresses two
questions posed by Patriksson \cite{Patriksson}.  First, it shows that the
sparsity can be exploited within the setting of a general-purpose optimizer.
Second, it provides an efficient method that also avoids the usual
assumptions (see \S\ref{section:descr}) imposed by pegging or breakpoint
search methods on the domain, monotonicity or strict convexity of $f_i$ and
$g_i$.

In the next section, we review the optimality conditions for
(\ref{obj})--(\ref{intervals}).  In \S\ref{sec:methods} we describe the
breakpoint search and interior point methods, along with details of their
implementation. Section \ref{sec:problems} lays out the problem instances
used for the computational tests, and the results are discussed in
\S\ref{sec:results}.

\section{Optimality conditions}
\label{section:descr}

In this study we make the following assumptions:
\begin{enumerate}
\item[A1.] The \emph{relaxed} problem, in which (\ref{gcon}) is replaced
    by $g(x)\leq b$, has no optimum with $g(x)<b$.
\item[A2.] The function $f_i$ is decreasing on $[l_i,u_i]$ and $g_i$ is
    increasing on $[l_i,u_i]$ with $g(l)<b<g(u)$.
\end{enumerate}
The randomly generated test instances of \S\ref{sec:problems} all satisfy
these assumptions, which are needed for breakpoint search but not for the
interior point method.

In practice, we are more interested in the relaxed problem mentioned in
Assumption A1. However, we can easily determine whether either assumption
holds if we know the intervals of monotonicity for each $f_i$ and $g_i$.
Indeed, many treatments of resource allocation problems include one or both
of these assumptions because they can be inexpensively enforced through some
combination of initialization, preprocessing, and data generation.

Assumptions A1--A2 imply that (\ref{obj})--(\ref{intervals}) and the relaxed
problem are equivalent and admit an optimal solution; they also guarantee
that the Slater constraint qualification holds for the relaxed problem. By
Lagrangian duality, necessary and sufficient optimality conditions for
(\ref{obj})--(\ref{intervals}) can therefore be expressed as follows:
$g(x)=b$ and, for some real number $\rho$, $x$ is a solution to the separable
optimization subproblem
\begin{equation}\label{brkptsubproblem}
  \text{minimize } f(x)+\rho g(x)
   \text{ subject to $l\leq x\leq u$.}
\end{equation}
The dual objective is
\begin{equation}\label{dualobj}
  \rho\mapsto-b\rho+\sum_{i=1}^n\min_{x_i\in [l_i,u_i]}[f_i(x_i)+\rho g_i(x_i)],
\end{equation}
which attains its maximum; moreover, any maximizer $\rho$ is necessarily
nonnegative. The subproblem (\ref{brkptsubproblem}) has coordinate-wise
optimality conditions given by
\begin{align*}
  f_i'(x_i)+\rho g_i'(x_i)=0,&\quad \text{ if } l_i<x_i<u_i,\\
  f_i'(x_i)+\rho g_i'(x_i)\geq0,&\quad\text{ if } x_i=l_i,\\
  f_i'(x_i)+\rho g_i'(x_i)\leq0,&\quad\text{ if } x_i=u_i.
\end{align*}
The left-hand sides give the Karush-Kuhn-Tucker multipliers for the bounds
$l_i\leq x_i$ and $x_i\leq u_i$, respectively, as
\begin{gather*}
\lambda_i:=\max\{0,-[f_i'(x_i)+\rho g_i'(x_i)]\},
\\
   \mu_i:=\max\{0,f_i'(x_i)+\rho g_i'(x_i)\}.
\end{gather*}
Letting $s:=u-x$ denote the vector of slack variables for the upper bounds on
$x$, we express the Karush-Kuhn-Tucker (KKT) conditions for
(\ref{obj})--(\ref{intervals}) as
\begin{gather}
  \nabla f(x)+\rho\nabla g(x)-\lambda+\mu=0,\label{KKT:gradient}\\
  x+s=u,\label{KKT:slack}\\
  x\geq l,\; \lambda\geq0,\; s\geq0,\; \mu\geq0,\label{KKT:nonneg}\\
  \diag(x-l)\lambda=0, \; \diag(s)\mu=0,\label{KKT:complement}\\
  g(x)=b.\label{KKT:gequation}
\end{gather}
Here $\diag(z)$ denotes the diagonal matrix whose diagonal entries are the
entries of the vector $z$.

The three solution frameworks discussed in \S\ref{sec:introduction} utilize
the optimality conditions in different ways:
\begin{itemize}
\item Pegging methods solve subproblems of the form
    (\ref{obj})--(\ref{gcon}), but for which some variables are held
    fixed while the bounds (\ref{intervals}) for all remaining variables
    are omitted.
\item Breakpoint search  methods maximize the dual objective
    (\ref{dualobj}) by solving a sequence of subproblems of the form
    (\ref{brkptsubproblem}) at various values of $\rho$.
\item Primal-dual interior point methods apply Newton's method to
    perturbations of the KKT system
    (\ref{KKT:gradient})--(\ref{KKT:gequation}).
\end{itemize}
The pegging and breakpoint search methods both benefit considerably when
minimization of $x_i\mapsto f_i(x_i)+\rho g_i(x_i)$ can be handled
efficiently.  Because we focus on problems for which breakpoint search
dominates pegging, we do not include pegging methods in this study. In fact,
the pegging approach is not even well-defined for some of the problems we
consider, because the pegging subproblems do not admit optimal solutions.

\section{Methods and implementation} \label{sec:methods}

In this section, we describe the two main approaches considered in our
computational study.

\subsection{Breakpoint search}\label{sec:brkpt}

Breakpoint search %(see Figure \ref{fig:brkptalg})
is based on the observation that the dual objective (\ref{dualobj}) is
concave and defined piecewise with a finite number of easily calculated
breakpoints. The derivative, or subdifferential, of this objective is
nonincreasing.  A binary search of the breakpoints therefore identifies
either one that is a root or a pair that most closely bracket a root.

There are at most $2n$ breakpoints, occurring at $\rho$-values where some
$x_i\mapsto f_i(x_i)+\rho g_i(x_i)$ attains its minimum over $[l_i,u_i]$ at
an endpoint $l_i$ or $u_i$. Equivalently, a breakpoint makes the derivative
$x_i\mapsto f'_i(x_i)+\rho g'_i(x_i)$ nonnegative at $l_i$ or nonpositive at
$u_i$.  Consequently, all breakpoints have the form
$\rho^+_i:=-f_i'(l_i)/g'_i(l_i)$ or $\rho^-_i:=-f_i'(u_i)/g'_i(u_i)$. The
monotonicity of $f_i$ and $g_i$ allow us to define $\rho^+_i=\infty$ when
$g'_i(l_i)=0$ and to guarantee that $g'_i(u_i)>0$ in the definition of
$\rho^-_i$. The convexity and monotonicity of $f_i$ and $g_i$ also guarantee
that $0\leq\rho^-_i\leq \rho^+_i$.

The binary search sequentially refines a bracketing $\rho^-<\rho^*<\rho^+$
until the true root $\rho^*$ lies between two consecutive breakpoints. The
bracket is adjusted inward by finding a breakpoint $\rho$ within it and
testing the sign of the derivative of the dual objective (\ref{dualobj}). To
evaluate that derivative at $\rho$, we first fix
\begin{equation} \label{xfixing}
x_i:=
 \begin{cases}
   l_i,&\text{if }\rho\geq\rho^+_i,\\
   u_i,&\text{if }\rho\leq\rho^-_i.
 \end{cases}
\end{equation}
The remaining minimizers are critical points: $f'_i(x_i)+\rho g'_i(x_i)=0$
and $l_i<x_i<u_i$.  Depending on the problem data, these critical points
might be found (a) in closed form, (b) by using a problem-specific
implementation of Newton's method, or (c) by means of a general-purpose
Newton's method with Armijo linesearch for sufficient decrease and damping
(as needed) to maintain $l_i<x_i<u_i$.  The derivative value at $\rho$ is
then given by $-b+\sum_ig_i(x_i)$, the sign of which determines whether
$\rho$ becomes the new $\rho^-$ or $\rho^+$. This in turn determines, through
(\ref{xfixing}), that some values of $x_i$ shall remain fixed and can
therefore be removed from further consideration.

The final bracket, if nontrivial, consists of two closest breakpoints with
the optimal value of $\rho$ lying somewhere between them. To interpolate
between them, our implementation finds $\rho$ and the unfixed
$x_i$-coordinates (denoted by $i \in I$) simultaneously by applying a
multi-dimensional Newton's method with Armijo linesearch to the corresponding
Lagrange multiplier conditions $\sum_{i\in I}g_i(x_i)=\hat b$ and
$f'_i(x_i)+\rho g'_i(x_i)=0$ for $i\in I$.

Throughout the procedure, the subproblem optimizations are initialized using
the corresponding solutions from prior iterations.  Also, we extract the
required median values without sorting the list of breakpoints in advance,
which can yield significant computational savings if each subproblem solution
requires only a few operations per index $i$
\cite{Brucker,Kiwiel08brkpt,Pardalos1990,Robinson}.

\subsection{Interior point method}\label{sec:ipm}

The primal-dual interior point method solves the KKT optimality conditions
(\ref{KKT:gradient})--(\ref{KKT:gequation}) for the variables
$(x,\lambda,s,\mu,\rho)$. Its operation preserves strict inequality for the
simple bounds (\ref{KKT:nonneg}), only allowing them to become active in the
limit. The method is based on the perturbed KKT system
\begin{gather}
  \nabla f(x)+\rho\nabla g(x)-\lambda+\mu=0,\label{IPM:gradient}\\
  x+s=u,\label{IPM:slack}\\
  \diag(x-l)\lambda=\tau e, \; \diag(s)\mu=\tau e,\label{IPM:complement}\\
  g(x)=b,\label{IPM:gequation}
\end{gather}
where $e$ denotes the $n$-vector of all ones and the inequalities $x>l$,
$\lambda>0$, $s>0$, $\mu>0$ are enforced separately. The system
(\ref{IPM:gradient})--(\ref{IPM:gequation}) is algebraically equivalent to
the Lagrange multiplier equations for a related optimization problem
involving barrier functions for the bounds (\ref{intervals}):
\begin{align*}
  &\text{minimize}\quad f(x)-\tau\sum_{i=1}^n[\ln(x_i-l_i)+\ln(s_i)]\\
  &\text{over all}\quad x>l,\;s>0\\
  &\text{subject to}\quad g(x)=b,\; x+s=u.
\end{align*}
As the barrier parameter $\tau>0$ is driven to zero, we expect the (unique)
solution $(x,\lambda,s,\mu,\rho)$ of
(\ref{IPM:gradient})--(\ref{IPM:gequation}) to tend toward the solution set
of the original KKT system (\ref{KKT:gradient})--(\ref{KKT:gequation}).

In each iteration of the interior point method, we calculate a Newton search
direction for the perturbed system
(\ref{IPM:gradient})--(\ref{IPM:gequation}) and then take a step along that
direction, damped so as to preserve $x>l$, $\lambda>0$, $s>0$, $\mu>0$. Next,
the value of $\tau$ is adjusted and the iteration repeats.  The method stops
when the residuals $r_d:=\nabla f(x)+\rho\nabla g(x)-\lambda+\mu$,
$r_l:=\diag(x-l)\lambda$, $r_u:=\diag(s)\mu$, and $r_g:=g(x)-b$ are small
enough.

Ours is a rudimentary implementation aimed at any nonlinear programming
formulation involving simple bounds and an equality constraint. The
algorithmic parameters were assigned values that gave reliable performance in
preliminary testing: all relative tolerances for residuals were set to
$10^{-10}$, $\tau$ was set to $0.25$ of the current duality gap
$(x\cdot\lambda+s\cdot\mu)/2n$, and the step size was taken to be the smaller
of unity or $0.8$ of the feasible step.  No attempt was made to provide
theoretical guarantees of global convergence, superlinear convergence, or
complexity.  However, the implementation correctly solved all the instances
described in \S\ref{sec:problems} and easily outperformed breakpoint search
on challenging problems of moderate to very large dimension. It therefore met
the needs of the present study. The key to making it competitive is the fact
that the linear system defining the Newton search direction can be solved in
$Cn$ arithmetic operations for a small fixed value of $C$, as we show next.

{\allowdisplaybreaks

To simplify the notation, we introduce a vector $h$ with entries
$h_i:=f_i''(x_i)+\rho g_i''(x_i)$ and let $\xi$ denote $x-l$. We also use
uppercase letters to denote these diagonal matrices: $\Xi:=\diag(\xi)$,
$\Lambda:=\diag(\lambda)$, $S:=\diag(s)$, $M:=\diag(\mu)$, $H:=\diag(h)$. The
linear system satisfied by the search direction is then
\begin{equation}
  \begin{bmatrix}
   H&-I&0&I&\nabla g\\
   \Lambda&\Xi&0&0&0\\
   0&0&M&S&0\\
   I&0&I&0&0\\
   \nabla g^T&0&0&0&0
   \end{bmatrix}
   \begin{bmatrix}
     \Delta x\\ \Delta\lambda \\ \Delta s\\ \Delta\mu \\ \Delta\rho
   \end{bmatrix}=
   \begin{bmatrix}
     r_d\\ r_l\\ r_u\\ 0\\ r_g
   \end{bmatrix}.
   \label{linearsystem}
\end{equation}
To solve (\ref{linearsystem}), first calculate the vectors
\begin{gather*}
w:=h+\Xi^{-1}\lambda+S^{-1}\mu,\\
y:=r_d+\Xi^{-1}r_l-S^{-1}r_u,\\
z:=W^{-1}\nabla g
\end{gather*}
and let $\eta:=-1/\nabla g^Tz$. Multiplying (\ref{linearsystem}) from the
left by
\[ \begin{bmatrix}
    W^{-1}&0&0&0&0\\
    0&I&0&0&0\\
    0&0&I&0&0\\
    0&0&0&I&0\\
    -\eta z^T&0&0&0&\eta
    \end{bmatrix}
    \begin{bmatrix}
    I&\Xi^{-1}&S^{-1}&S^{-1}M&0\\
    0&\Xi^{-1}&0&0&0\\
    0&0&S^{-1}&0&0\\
    0&0&0&I&0\\
    0&0&0&0&1
    \end{bmatrix}
\]
yields
\[ \begin{bmatrix}
   I&0&0&0&z\\
   \Xi^{-1}\Lambda&I&0&0&0\\
   0&0&S^{-1}M&I&0\\
   I&0&I&0&0\\
   0&0&0&0&1
   \end{bmatrix}
   \begin{bmatrix}
     \Delta x\\ \Delta\lambda \\ \Delta s\\ \Delta\mu \\ \Delta\rho
   \end{bmatrix}=
   \begin{bmatrix}
     W^{-1}y\\ \Xi^{-1}r_l\\ S^{-1}r_u\\ 0\\ \eta(r_g-z^Ty)
   \end{bmatrix}
\]
from which we can read off the solution to (\ref{linearsystem}) as:
\begin{gather*}
\Delta\rho =\eta(r_g-z^Ty),\\
\Delta x=W^{-1}y-(\Delta\rho)z,\\
\Delta s=-\Delta x,\\
\Delta\lambda=\Xi^{-1}r_l-\Xi^{-1}\Lambda\Delta x,\\
\Delta\mu=S^{-1}r_u-S^{-1}M\Delta s.
\end{gather*}
The saved values of $\Xi^{-1}r_l$, $\Xi^{-1}\Lambda$, $S^{-1}r_u$, $S^{-1}M$
from the calculation of $w$ and $y$ can be reused here. The solution of
(\ref{linearsystem}) requires $9n-1$ additions/subtractions, $5n+1$
multiplications, and $6n+1$ divisions.

Table \ref{tab:linsys} shows that the proposed method for solving
(\ref{linearsystem}) is much faster than a standard linear solver, namely,
the MATLAB sparse LU factorization with approximate minimum-degree reordering
of columns.  The proportionate speed-up seen here completely accounts for the
superiority of our interior point method over the general breakpoint search
(see \S\ref{sec:results}).
\begin{table}[!tbh]\centering
\caption{Solution time in ms for linear system (\ref{linearsystem}).}
\label{tab:linsys}
\begin{tabular}{lrrrrrr}
\hline
%&\hbox to0pt{\hss$n={}$}$10^2$& $10^3$& $10^4$& $10^5$& $10^6$\\
&\multicolumn{5}{c}{dimension}\\
\cline{2-6}
solver&$10^2$& $10^3$& $10^4$& $10^5$& $10^6$\vrule width0pt height10pt\\
\hline
proposed &  0.03 &  0.08 &   0.50  &    5.53 &    69 \\
LU colamd &  0.85 &  3.25 &  43.16 &  526.98 &  5627 \\
%proposed &  0.03 &  0.08 &   0.50 &    5.53 &    69.37 \\
%LU colamd &  0.85 &  3.25 &  43.16 &  526.98 &  5626.91 \\
\hline
\end{tabular}
\end{table}

}

\section{Test instances} \label{sec:problems}

For our computational tests, we selected five problem classes involving
mathematical forms of potential interest in operations research.  None admits
closed-form solutions for its separable breakpoint subproblems. Instances
were generated so that assumptions A1--A2 of \S\ref{section:descr} were
satisfied, after possible reorientation of intervals. In the following, the
notation $z\sim U(a,b)$ indicates that the value $z$ was selected according
to a continuous uniform distribution on the open interval $(a,b)$, whereas
$z\sim N(\mu,\sigma)$ indicates that $z$ was selected according to a normal
distribution with mean $\mu$ and standard deviation $\sigma$.

\subsection{Resource renewal}

Problems in this class have $f_i(x_i)=a_ix_i(e^{-1/x_i}-1)$ and
$g_i(x_i)=c_ix_i$ for $x_i>0$, as studied by Melman and Rabinowitz
\cite{Melman}. For convenience, we extend $f_i$ to a $C^\infty$ convex
function on the real line by defining $f_i(x_i)=-x_i$ for $x_i\leq0$.
Instances were generated as follows:
\begin{itemize}
\item $a_i,c_i\sim U(0.001,1000)$;
\item $b=1.1\sum_ic_i\xi_i$, where $\gamma=\min_j\{a_j/c_j\}$ and
\[\xi_i=\begin{cases}
    0,&\text{if $a_i/c_i>\gamma$,}\\
    \mathop{\rm argmin}\limits_{x_i} f_i(x_i)+\gamma g_i(x_i), &\text{if
      $a_i/c_i\leq\gamma$;}
    \end{cases}
\]
\item $l_i=0$ and $u_i=b/c_i$.
\end{itemize}

\subsection{Weighted $p$-norm over a ball}
\label{sec:problems:pnormrball}

Problems in this class have $f_i(x_i)=a_i|x_i-y_i|^p$ and $g_i(x_i)=|x_i|^r$.
Note that $f_i$ and $g_i$ are everywhere twice differentiable when
$p,r\geq2$. Instances with $p,r\in\{2,2.5,3,4\}$ and $p\not=r$ were generated
as follows:
\begin{itemize}
\item $a_i\sim U(1,10)$;
\item $l_i\sim U(0,5)$, $u_i\sim U(l_i,l_i+5)$;
\item $y_i\sim U(u_i,u_i+5)$, $b\sim U(g(l),g(u))$.
\end{itemize}

\subsection{Sums of powers}

Problems in this class have $f_i(x_i)=a_i|x_i-y_i|^{p_i}$ and
$g_i(x_i)=|x_i|^{r_i}$.  Instances were generated as follows:
\begin{itemize}
\item $a_i\sim U(1,10)$;
\item $p_i,r_i\sim U(2,4)$;
\item $l_i\sim U(0,5)$, $u_i\sim U(l_i,l_i+5)$;
\item $y_i\sim U(u_i,u_i+5)$, $b\sim U(g(l),g(u))$.
\end{itemize}

\subsection{Convex quartic over a simplex}

This class of problems has $f_i(x_i)=a_ix_i^4+b_ix_i^3+c_ix_i^2+d_ix_i$ and
$g_i(x_i)=x_i$.  Instances of these problems were generated as follows:
\begin{itemize}
\item $a_i=(\xi_i^2+\eta_i^2)/\sqrt8$,
    $b_i=(\xi_i\zeta_i+\eta_i\chi_i)/\sqrt3$ and
    $c_i=(\zeta_i^2+\chi_i^2)/\sqrt8$, with
    $\xi_i,\eta_i,\zeta_i,\chi_i\sim N(0,1)$;
\item $d_i=-f'_i(\tau_i)$, with $\tau_i\sim U(0,10)$;
\item $u_i=\min(\tau_i,\lambda_i)$, with $\lambda_i\sim U(0,\tau_i)$;
\item $l_i\sim U(0,u_i)$;
\item $b\sim U(g(l),g(u))$.
\end{itemize}
The choice of coefficients for $f_i$ guarantees that $f_i$ is strictly convex
on the real line, which is true if and only if $8a_ic_i>3b_i^2$, $a_i>0$, and
$c_i>0$.  Equivalently, the matrix
\[ \left[\begin{matrix}
    \sqrt{8}a_i&\sqrt{3}b_i\\
    \sqrt{3}b_i&\sqrt{8}c_i
   \end{matrix}\right]
\]
must be positive definite.  This can be ensured by selecting $a_i,b_i,c_i$ to
be the rescaled entries of a matrix formed as $A^TA$, where the entries of
$A$ are given by $\xi_i,\eta_i,\zeta_i,\chi_i$.  The choice of $d_i$
guarantees that the critical point of $f_i$ is positive, after which the
bounds are chosen so that each $f_i$ has the same monotonicity.

\subsection{Log-exponential}

Problems in this class have $g_i(x_i)=c_ix_i$ and
\[f_i(x_i)=\ln\left[\sum_{j=1}^5\exp(a_{ij}x_i+d_{ij})\right].\]
Instances were generated as follows:
\begin{itemize}
\item $d_{ij}\sim N(0,1)$ and $c_i\sim U(0,10)$;
\item $\xi_{ij}\sim N(0,1)$ and $\displaystyle\zeta_i\sim \begin{cases}
    U(0,1),&\text{if $i\in I$},\\ N(0,1),&\text{if $i\not\in I$};
    \end{cases}$
\item $\displaystyle a_{ij}=
   \begin{cases}
    |\xi_{ij}|,&\text{if $\xi_{ij}>0$ for all $j$,}\\
    \xi_{ij},&\text{otherwise;}
   \end{cases}
   $
\item $\chi_i=\mathop{\rm argmin}f_i$;
\item $\displaystyle u_i=
    \begin{cases}
      \min(\chi_i,
         1.2\zeta_i\chi_i),&\text{if $i\in I$},\\
      5\zeta_i,&\text{if $i\not\in I$};
    \end{cases}$
\item $l_i=u_i-0.05|u_i|-5|\eta_i|$, with $\eta_i\sim N(0,1)$;
\item $b\sim U(g(l),g(u))$.
\end{itemize}

\subsection{A note on easier problems}

In addition to the problems described above, we ran similar tests on randomly
generated instances from several classes of problems admitting closed-form
algebraic solutions to the either pegging or breakpoint subproblems.  We
don't report the results on these easier problems beyond the following brief
summary. Unsurprisingly, problem-specific pegging methods (when available)
handily outperformed all other approaches. However, the results were mixed
concerning problem-specific breakpoint search versus the interior point
method: the former was faster for low- to medium-dimensional instances but
lost its edge as the dimension increased, so that the interior point method
was generally faster for $n\geq10^4$ or $n\geq10^5$. Regardless of the
problem class, the performance of the interior point method relative to the
general breakpoint search of \S\ref{sec:brkpt} was similar to the results
described in \S\ref{sec:results} below.

\section{Computational results} \label{sec:results}

The procedures of \S\ref{sec:methods} were coded in MATLAB 7.12 and their
performance compared on randomly generated instances as described in
\S\ref{sec:problems}.  All tests were performed on a dedicated
2$\times$quad-core Intel 64-bit (2.26 GHz) platform with 24GB RAM running
CentOS Linux.

We attempted to solve instances of dimension $10^k$ for $k\in\{2,3,4,5,6\}$
with both methods for each problem class.  One hundred random instances were
generated at each dimension for four of the five problem classes. The
exception was the class described in \S\ref{sec:problems:pnormrball}, for
which we generated 100 instances at each dimension for each value of
$p,r\in\{2,2.5,3,4\}$. Performance differences among these combinations of
$p$ and $r$ were detectable, but too small to warrant separating out the
results for discussion.  We therefore report them in aggregate over all $p$
and $r$, and remark that the higher values of $p$ or $r$ tended to require
slightly longer running times than did the smaller values.

Based on preliminary testing, an \emph{a priori} time limit of $10^{k-2}$
seconds was imposed on each attempt at solution. The interior point method
was run first and never exceeded this time limit.  Consequently, for higher
dimensional instances on some problem types, the breakpoint searches were
limited to at most a factor of ten over the worst runtime for the interior
point method on problems of the same type and size.

\begin{figure*}[!tbh]\centering
\includegraphics[scale=.80,draft=false]{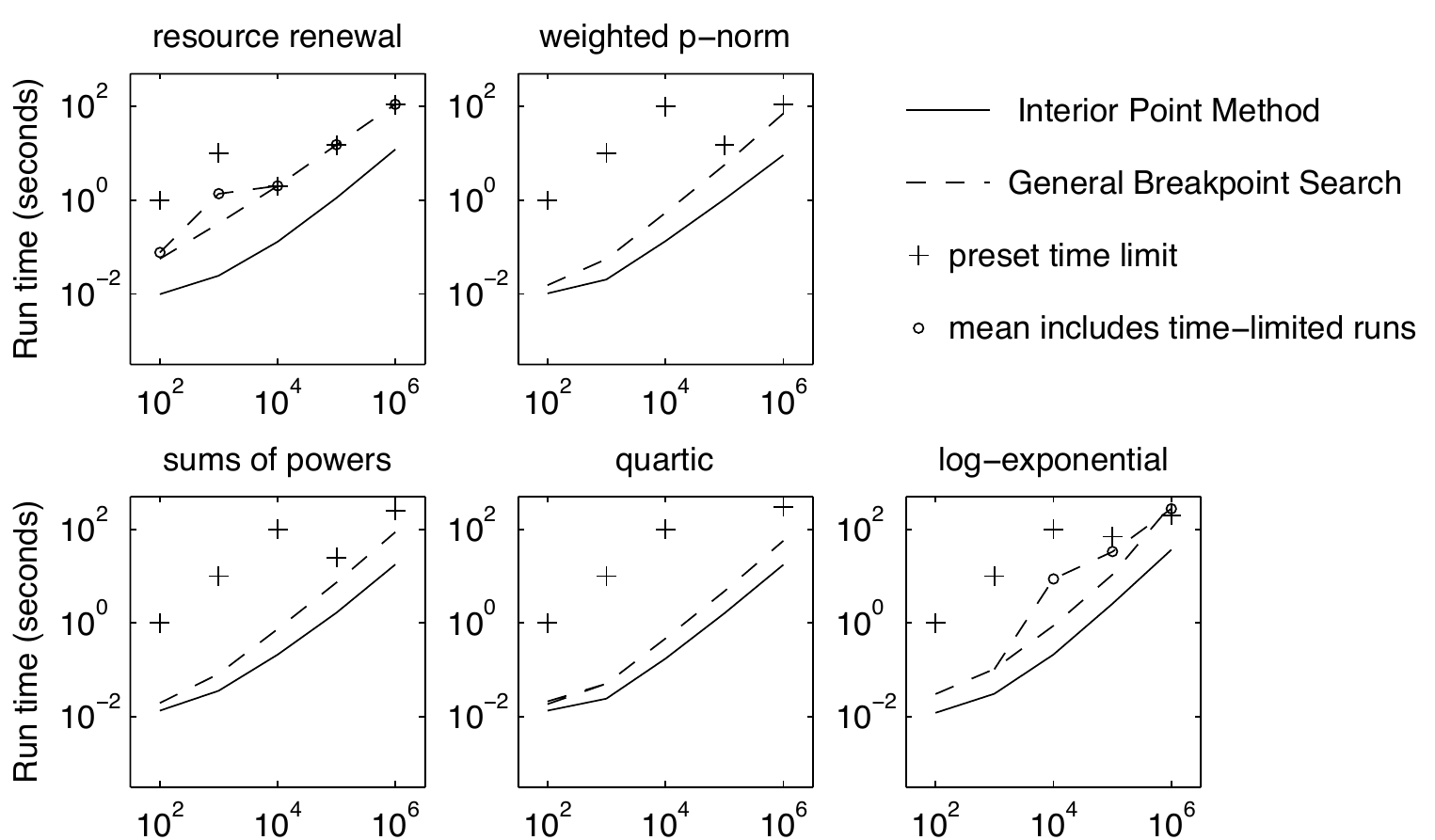} % arxiv fix
\caption{Mean and median running times (seconds) for three methods on ten
problem classes.} \label{fig:runtimes}
\end{figure*}

The results of the tests are presented by problem class in Figure
\ref{fig:runtimes}.
Because the breakpoint search often exceeded the given time limits, the
graphs consist mainly of the median runtimes.  Mean runtimes are drawn
separately whenever they can be visually distinguished from the medians on
the scale shown; the mean curves are always the upper branches when two
curves of the same line type are shown.  Means that include runtimes at their
limits are specially marked.

The interior point method clearly dominates the general breakpoint search,
often by an order of magnitude. Table \ref{tab:wins} shows the frequency with
which such dominance occurs. The common scale and position on the five graphs
in Figure \ref{fig:runtimes} suggest that running times for the interior
point method do not depend greatly on the specific type problem (aside from
the expense of function evaluations). On the other hand, the performance gap
between the two optimization methods is smaller than the gap between the two
linear-system solvers considered in \S\ref{sec:ipm}.  We conclude that when
algebraic simplifications due to the form of $f_i$ and $g_i$ are unavailable,
an interior point method is a strong option for solving problems of the form
(\ref{obj})--(\ref{intervals}). However, its competitiveness relies even more
heavily than usual on the efficiency of the underlying linear-system solver.

\begin{table}[!tbh]\centering
\caption{Win percentage of IPM over general breakpoint
method}\label{tab:wins}
\begin{tabular}{lrrrrrr}
\hline
%problem&\hbox to0pt{\hss$n={}$}$10^2$& $10^3$& $10^4$& $10^5$& $10^6$\\
&\multicolumn{5}{c}{dimension}\\
\cline{2-6}
problem class&$10^2$& $10^3$& $10^4$& $10^5$& $10^6$\vrule width0pt height10pt\\
\hline
     resource renewal  &   0 &  73 &  86 & 100 & 100\\
     weighted $p$-norm &   4 &  95 &  99 & 100 & 100\\
     sums of powers    &   4 &  97 & 100 & 100 & 100\\
     quartic           &  86 & 100 & 100 & 100 &  97\\
     log-exponential   & 100 & 100 & 100 & 100 & 100\\
\hline
\end{tabular}
\end{table}

\section*{Acknowledgements}

The authors thank Professors A.~John Bailer and Douglas A.~Noe for
suggestions on improving the paper. Miami University's ``RedHawk'' computing
cluster was used for all computational tests reported here.

\end{document}